# Scale space consistency of piecewise constant least squares estimators – another look at the regressogram

**Leif Boysen** [1,*], **Volkmar Liebscher**[2,†], **Axel Munk**[1] **and Olaf Wittich**[3,†]

*Universität Göttingen, Universität Greifswald, Universität Göttingen and Technische Universiteit Eindhoven*

**Abstract:** We study the asymptotic behavior of piecewise constant least squares regression estimates, when the number of partitions of the estimate is penalized. We show that the estimator is consistent in the relevant metric if the signal is in $L^2([0,1])$, the space of *càdlàg* functions equipped with the Skorokhod metric or $C([0,1])$ equipped with the supremum metric. Moreover, we consider the family of estimates under a varying smoothing parameter, also called scale space. We prove convergence of the empirical scale space towards its deterministic target.

## 1. Introduction

Initially, the use of piecewise constant functions for regression has been proposed by [25], who called the corresponding reconstruction the regressogram. [25] proposed it as a simple exploratory tool. For a given set of jump locations, the regressogram simply averages the data between two successive jumps. A difficult issue, however, is a proper selection of the location of jumps and its convergence analysis.

Approximation by step functions is well examined in approximation theory (see e.g., [7]), and there are several statistical estimation procedures which use locally constant reconstructions. [14] studied the case where the signal is a step function with one jump and showed that in this case the signal can be estimated at the parametric $n^{-1/2}$-rate and that the jump location can be estimated at a rate of $n^{-1}$. This was generalized by [28] and [29] to step functions with a given a known upper bound for the number of jumps. The locally adaptive regression splines method by [16] and the taut string procedure by [6] use locally constant estimates to reconstruct unknown regression functions, which belong to more general function classes. Both methods reduce the complexity of the reconstruction by minimizing the total variation of the estimator, which in turn leads to a small number of local extreme values.

*Supported by Georg Lichtenberg program "Applied Statistics & Empirical Methods" and DFG graduate program 1023 "Identification in Mathematical Models".

†Supported in part by DFG, Sonderforschungsbereich 386 "Statistical Analysis of Discrete Structures".

‡Supported by DFG grant "Statistical Inverse Problems under Qualitative Shape Constraints".
[1]Institute for Mathematical Stochastics, Georgia Augusta University Goettingen, Maschmuehlenweg 8-10, D-37073 Goettingen, Germany, e-mail: stochastik@math.uni- goettingen.de; boysen@math.uni-goettingen.de
[2]University Greifswald.
[3]Technical University Eindhoven.

In this work we choose a different approach and define the complexity of the reconstruction by the number of intervals where the reconstruction is constant, or equivalently by the number of jumps of the reconstruction. Compared to the total variation approach, this method obviously captures extreme plateaus more easily but is less robust to outliers. This might be of interest in applications where extreme plateaus are informative, like for example in mass spectroscopy.

Throughout the following, we assume a regression model of the type

$$(1) \qquad Y_{i,n} = \bar{f}_{i,n} + \xi_{i,n}, \qquad (i = 1, \ldots, n),$$

where $(\xi_{i,n})_{i=1,\ldots,n}$ is a triangular array of independent zero-mean random variables and $\bar{f}_{i,n}$ is the mean value of a square integrable function $f \in L^2([0,1))$ over the interval $[(i-1)/n, i/n]$ (see e.g. [9]),

$$(2) \qquad \bar{f}_{i,n} = n \int_{(i-1)/n}^{i/n} f(u)\, du.$$

This model is well suited for physical applications, where observations of this type are quite common.

We consider minimizers $T_\gamma(Y_n) \in \operatorname{argmin} H_\gamma(\cdot, Y_n)$ of the hard thresholding functional

$$(3) \qquad H_\gamma(u, Y_n) = \gamma \cdot \#J(u) + \frac{1}{n}\sum_{i=1}^{n}(u_i - Y_{i,n})^2,$$

where

$$J(u) = \{i : 1 \le i \le n-1, u_i \ne u_{i+1}\}$$

is the set of *jumps* of $u$. In the following we will call the minimizers of (3) jump penalized least squares estimators or short *Jplse*.

Clearly choosing $\gamma$ is equivalent to choosing a number of partitions of the *Jplse*. Figure 1 shows the *Jplse* for a sample dataset and different choices of the smoothing parameter $\gamma$.

This paper complements work of the authors on convergence rates of the *Jplse*. [2] show that given a proper choice of the smoothing parameter $\gamma$ it is possible to obtain optimal rates for certain classes of approximation spaces under the assumption of subgaussian tails of the error distribution. As special cases the class of piecewise Hölder continuous functions of order $0 < \alpha \le 1$ and the class of functions with bounded total variation are obtained.

In this paper we show consistency of regressograms constructed by minimizing (3) for arbitrary $L^2$ functions and more general assumptions on the error. If the true function is *càdlàg*, we additionally show consistency in the Skorokhod topology. This is a substantially stronger statement than the $L_2$ convergence and yields consistency of the whole *graph* of the estimator.

In concrete applications the choice of the regularization parameter $\gamma > 0$ in (3), which controls the degree of smoothness (which means just the number of jumps) of the estimate $T_\gamma(Y_n)$, is a delicate and important task. As in kernel regression [18, 23], a screening of the estimates over a larger region can be useful (see [16, 26]). Adapting a viewpoint from computer vision (see [15]), [3, 4] and [17] proposed to consider the family $(T_\gamma(f))_{\gamma>0}$, denoted as *scale space*, as target of inference. This was justified in [4] by the fact that the empirical scale space converges towards that of the actual density or regression function pointwisely and uniformly on compact sets. The main motivation for analyzing the scale space is exploration of structures as peaks and valleys in regression and detection of modes in density estimation. Properties of the scale space in kernel smoothing are that structures like modes disappear monotonically for a shrinking resolution level and that the reconstruction changes continuously with respect to the bandwidth. For the *Jplse*, the family $(T_\gamma(f))_{\gamma>0}$ behaves quite differently. Notable distinctions are that jumps may not change monotonically and that there are only finitely many possible different estimates. To deal with these features, we consider convergence of the scale space in the space of *càdlàg* functions equipped with the Skorokhod $J_1$ topology. In this setting we deduce (under identifiability assumptions) convergence of the empirical scale space towards its deterministic target. Note that the computation of the empirical scale space is feasible. The family $(T_\gamma(Y_n))_{\gamma>0}$ can be computed in $O(n^3)$ and the minimizer for one $\gamma$ in $O(n^2)$ steps (see [26]).

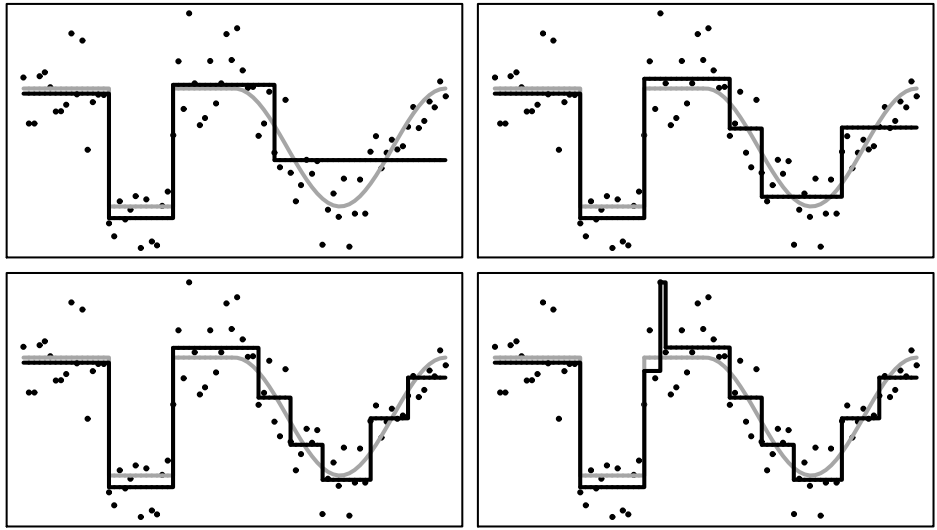

FIG 1. *The* Jplse *for different values of $\gamma$. The dots represent the noisy observations of some signal $f$ represented by the grey line. The black line shows the estimator, with $\gamma$ chosen such that the reconstruction has four, six, eight and ten partitions, respectively.*

The paper is organized as follows. After introducing some notation in Section 2, we provide in Section 3.1 the consistency results for general functions in the $L^2$ metric. In Section 3.2 we present the results of convergence in the Skorokhod topology. Finally in Section 3.3 convergence results for the scale space are given. The proofs as well as a short introduction to the concept of epi-convergence, which is required in the main part of the proofs, are given in the Appendix.

## 2. Model assumptions

By $S([0,1)) = \text{span}\{\mathbf{1}_{[s,t)} : 0 \leq s < t \leq 1\}$ we will denote the space of step functions with a finite but arbitrary number of jumps and by $D([0,1))$ the *càdlàg* space of right continuous functions on $[0,1]$ with left limits and left continuous at 1. Both will be considered as subspaces of $L^2([0,1))$ with the obvious identification of a function with its equivalence class, which is injective for these two spaces. More generally, by $D([0,1),\Theta)$ and $D([0,\infty),\Theta)$ we will denote spaces of functions with values in a metric space $(\Theta, \rho)$, which are right continuous and have left limits. $\|\cdot\|$ will denote the norm of $L^2([0,1))$ and the norm on $L^\infty([0,1))$ is denoted by $\|\cdot\|_\infty$.

Minimizers of the hard thresholding functionals (3) will be embedded into $L^2([0,1))$ by the map $\iota^n : \mathbb{R}^n \longmapsto L^2([0,1))$,

$$\iota^n((u_1, \ldots, u_n)) = \sum_{i=1}^n u_i \mathbf{1}_{[(i-1)/n, i/n)}.$$

Under the regression model (1), this leads to estimates $\hat{f}_n = \iota^n(T_{\gamma_n}(Y_n))$, i.e.

$$\hat{f}_n \in \iota^n(\operatorname{argmin} H_{\gamma_n}(\cdot, Y_n)).$$

Note that, for a functional $F$ we denote by $\operatorname{argmin} F$ the whole set of minimizers. Here and in the following $(\gamma_n)_{n \in \mathbb{N}}$ is a (possibly random) sequence of smoothing parameters. We suppress the dependence of $\hat{f}_n$ on $\gamma_n$ since this choice will be clear from the context.

For the noise, we assume the following condition.

(A) *For all $n \in \mathbb{N}$ the random variables $(\xi_{i,n})_{1 \leq i \leq n}$ are independent. Moreover, there exists a sequence $(\beta_n)_{n \in \mathbb{N}}$ with $n^{-1}\beta_n \to 0$ such that*

(4) $$\max_{1 \leq i \leq j \leq n} \frac{(\xi_{i,n} + \cdots + \xi_{j,n})^2}{j - i + 1} \leq \beta_n \quad \mathbb{P}\text{-}a.s.,$$

*for almost every $n$.*

The behavior of the process (4) is well known for certain classes of i.i.d. subgaussian random variables (see e.g. [22]). If for example $\xi_{i,n} = \xi_i \sim N(0, \sigma^2)$ for all $i = 1, \ldots, n$ and all $n$, we can choose $\beta_n = 2\sigma^2 \log n$ in Condition (A). The next result shows that (A) is satisfied for a broad class of subgaussian random variables.

**Lemma 1.** *Assume the noise satisfies the following generalized subgaussian condition*

(5) $$\mathbb{E}e^{\nu \xi_{i,n}} \leq e^{\alpha n^\zeta \nu^2}, \qquad (\text{for all } \nu \in \mathbb{R}, n \in \mathbb{N}, 1 \leq i \leq n)$$

*with $0 \leq \zeta < 1$ and $\alpha > 0$. Then there exist a $C > 0$ such that for $\beta_n = Cn^\zeta \log n$ Condition (A) is satisfied.*

A more common moment condition is given by the following lemma.

**Lemma 2.** *Assume the noise satisfies*

(6) $$\sup_{i,n} \mathbb{E}|\xi_{i,n}|^{2m} < \infty, \qquad (\text{for all } n \in \mathbb{N}, 1 \leq i \leq n)$$

*for $m > 2$. Then for all $C > 0$ and $\beta_n = C(n \log n)^{2/m}$ Condition (A) is satisfied.*

## 3. Consistency

In order to extend the functional in (3) to $L^2([0,1))$, we define for $\gamma > 0$, the functionals $H_\gamma^\infty : L^2([0,1)) \times L^2([0,1)) \longmapsto \mathbb{R} \cup \infty$:

$$H_\gamma^\infty(g, f) = \begin{cases} \gamma \cdot \#\mathcal{J}(g) + \|f - g\|^2, & g \in S([0,1)), \\ \infty, & \text{otherwise.} \end{cases}$$

Here
$$\mathcal{J}(g) = \{t \in (0,1) : g(t-) \neq g(t+)\}$$
is the set of jumps of $g \in S([0,1))$. For $\gamma = 0$, we set $H_0^\infty(g,f) = \|f - g\|^2$ for all $g \in L^2([0,1))$. The following lemma guarantees the existence of a minimizer.

**Lemma 3.** *For any $f \in L^2([0,1))$ and all $\gamma \geq 0$ we have*
$$\operatorname{argmin} H_\gamma^\infty(\cdot, f) \neq \varnothing.$$

In the following we assume that $Y_n$ is determined through (1), the noise $\xi_n$ satisfies (A) and $(\beta_n)_{n \in \mathbb{N}}$ is a sequence with $\beta_n/n \to 0$ such that (4) holds.

### 3.1. Convergence in $L^2$

We start with investigating the asymptotic behavior of the *Jplse* when the sequence $\gamma_n$ converges to a constant $\gamma$ greater than zero. In this case we do not recover the original function in the limit, but a parsimonious representation at a certain scale of interest determined by $\gamma$.

**Theorem 1.** *Suppose that $f \in L^2([0,1))$ and $\gamma > 0$ are such that $f_\gamma$ is a unique minimizer of $H_\gamma^\infty(\cdot, f)$. Then for any (random) sequence $(\gamma_n)_{n \in \mathbb{N}} \subset (0, \infty)$ with $\gamma_n \to \gamma$ $\mathbb{P}$-a.s., we have*
$$\hat{f}_n \xrightarrow[n \to \infty]{L^2([0,1))} f_\gamma \quad \mathbb{P}\text{-a.s.}$$

The next theorem states the consistency of the *Jplse* towards the true signal for $\gamma = 0$ under some conditions on the sequence $\gamma_n$.

(H)  $(\gamma_n)_{n \in \mathbb{N}}$ *satisfies* $\gamma_n \to 0$ *and* $\gamma_n n/\beta_n \to \infty$ $\mathbb{P}$-a.s..

**Theorem 2.** *Assume $f \in L^2([0,1))$ and $(\gamma_n)_{n \in \mathbb{N}}$ satisfies (H). Then*
$$\hat{f}_n \xrightarrow[n \to \infty]{L^2([0,1))} f, \quad \mathbb{P}\text{-a.s.}$$

### 3.2. Convergence in Skorokhod topology

As we use *càdlàg* functions for reconstructing the original signal, it is natural to ask, whether it is possible to obtain consistency in the Skorokhod topology.

We remember the definition of the Skorokhod metric [12, Section 5 and 6]. Let $\Lambda_\infty$ denote the set of all strictly increasing continuous functions $\lambda : \mathbb{R}_+ \mapsto \mathbb{R}_+$ which are onto. We define for $f, g \in D([0, \infty), \Theta)$
$$\rho(f(\lambda(t) \wedge u), g(t))$$
where $\mathrm{L}(\lambda) = \sup_{s \neq t \geq 0} |\log \frac{\lambda(t) - \lambda(s)}{t-s}|$. Similarly, $\Lambda_1$ is the set of all strictly increasing continuous onto functions $\lambda : [0,1] \mapsto [0,1]$ with appropriate definition of L. Slightly abusing notation, we set for $f, g \in D([0,1), \Theta)$,
$$\rho_S(f, g) = \inf \left\{ \max(\mathrm{L}(\lambda), \sup_{0 \leq t \leq 1} \rho(f(\lambda(t)), g(t))) : \lambda \in \Lambda_1 \right\}.$$

The topology induced by this metric is called $J_1$ topology. After determining the metric we want to use, we find that in the situation of Theorem 1 we can establish consistency without further assumptions, whereas in the situation of Theorem 2 $f$ has to belong to $D([0,1))$.

**Theorem 3.**  (i) *Under the assumptions of Theorem 1,*

$$\hat{f}_n \xrightarrow[n\to\infty]{D([0,1))} f_\gamma \quad \mathbb{P}\text{-}a.s.$$

(ii) *If $f \in D([0,1))$ and $(\gamma_n)_{n\in\mathbb{N}}$ satisfies (H), then*

$$\hat{f}_n \xrightarrow[n\to\infty]{D([0,1))} f \quad \mathbb{P}\text{-}a.s.$$

*If $f$ is continuous on $[0,1]$, then*

$$\hat{f}_n \xrightarrow[n\to\infty]{L^\infty([0,1])} f \quad \mathbb{P}\text{-}a.s.$$

### 3.3. Convergence of the scale spaces

As mentioned in the introduction, following [4], we now want to study the scale space family $(T_\gamma(f))_{\gamma>0}$ as target for inference. First we show that the map $\gamma \mapsto T_\gamma(f)$ can be chosen piecewise constant with finitely many jumps.

**Lemma 4.** *Let $f \in L^2([0,1))$. Then there exists a number $m(f) \in \mathbb{N} \cup \{\infty\}$ and a decreasing sequence $(\gamma_m)_{m=0}^{m(f)} \subset \mathbb{R} \cup \infty$ such that*

(i) $\gamma_0 = \infty, \gamma_{m(f)} = 0$,
(ii) *for all $1 \leq i \leq m(f)$ and $\gamma', \gamma'' \in (\gamma_i, \gamma_{i-1})$ we have that*

$$\operatorname{argmin} H^\infty_{\gamma'}(\cdot, f) = \operatorname{argmin} H^\infty_{\gamma''}(\cdot, f),$$

(iii) *for all $1 \leq i \leq m(f) - 1$ and $\gamma_{i+1} < \gamma' < \gamma_i < \gamma'' < \gamma_{i-1}$ we have:*

$$\operatorname{argmin} H^\infty_{\gamma_i}(\cdot, f) \supseteq \operatorname{argmin} H^\infty_{\gamma'}(\cdot, f) \cup \operatorname{argmin} H^\infty_{\gamma''}(\cdot, f),$$

*and*
(iv) *for all $\gamma' > \gamma_1$*

$$\operatorname{argmin} H^\infty_\infty(\cdot, f) = \operatorname{argmin} H^\infty_{\gamma'}(\cdot, f) = \{T_\infty(f)\}.$$

*Here $T_\infty(f)$ is defined by $T_\infty(f)(x) = \int f(u)\, du \mathbf{1}_{[0,1)}(x)$.*

Thus we may consider functions $\hat{\tau}_n \in D([0,\infty), L^2([0,1)))$ with

$$\hat{\tau}_n(\zeta) \in \iota^n(\operatorname{argmin} H_{1/\zeta}(\cdot, Y_n)),$$

for all $\zeta \geq 0$. We will call $\hat{\tau}_n$ the empirical scale space. Similarly, we define the deterministic scale space $\tau$ for a given function $f$, such that

(7) $$\tau(\zeta) \in \operatorname{argmin} H^\infty_{1/\zeta}(\cdot, f)), \qquad (\text{for all } \zeta \geq 0).$$

The following theorem shows that the empirical scale space converges almost surely to the deterministic scale space. Table 1 and Figure 2 demonstrate this in a finite setting for the blocks signal, introduced by [10].

**Theorem 4.** *Suppose $f \in L^2([0,1))$ is such that $\#\operatorname{argmin} H^\infty_\gamma(\cdot, f) = 1$ for all but a countable number of $\gamma > 0$ and $\#\operatorname{argmin} H^\infty_\gamma(\cdot, f) \leq 2$ for all $\gamma \geq 0$. Then $\tau$ is uniquely determined by (7). Moreover,*

$$\hat{\tau}_n \xrightarrow[n\to\infty]{} \tau \quad \mathbb{P}\text{-}a.s.$$

*holds both in $D([0,\infty), D([0,1)))$ and $D([0,\infty), L^2([0,1)))$.*

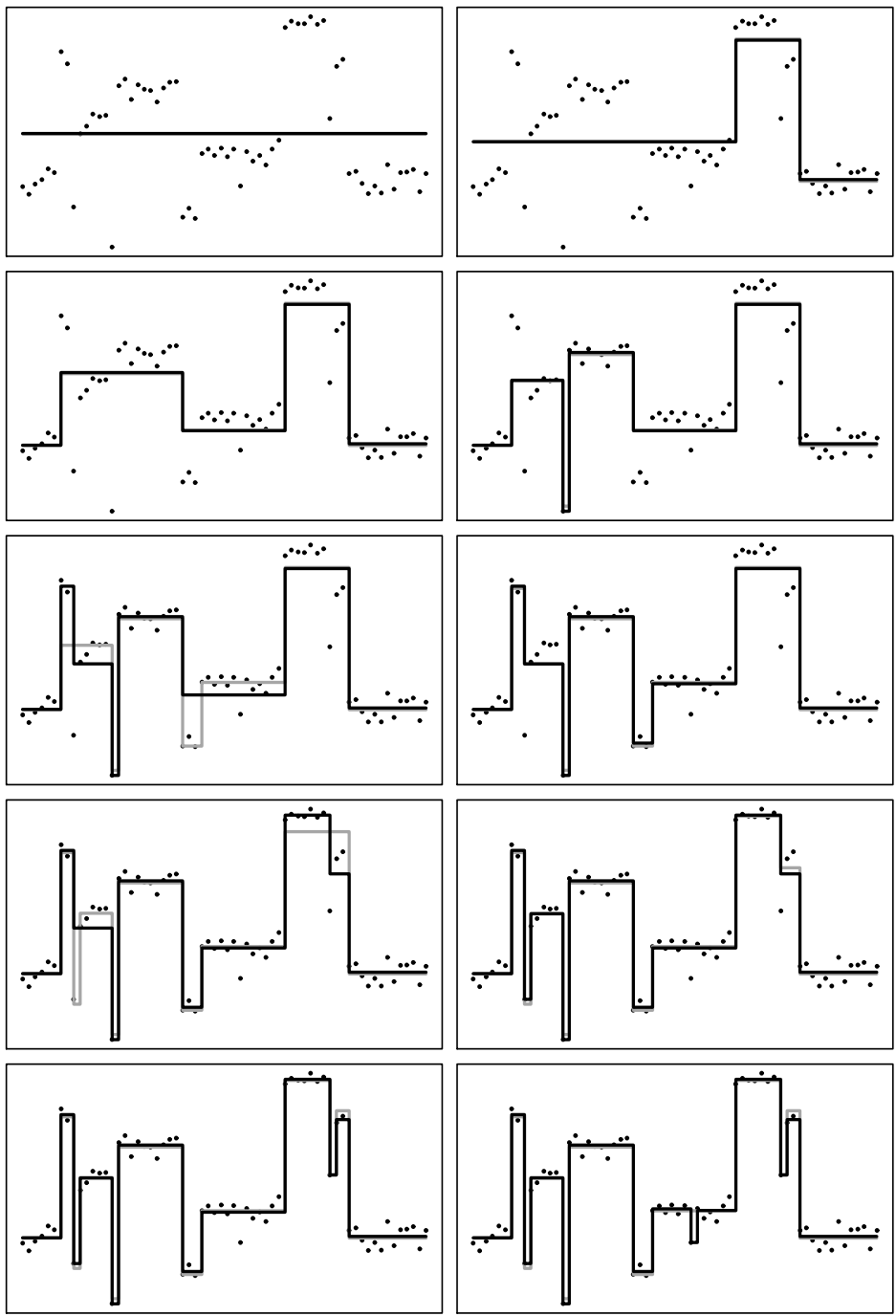

FIG 2. *Comparison of scale spaces. The "Blocks" data of [11] sampled at 64 points (dots) are compared with the different parts of the scale space derived both from the data (black) and the original signal (grey), starting with $\gamma = \infty$ and lowering its value from left to right and top to bottom. Note that for the original sampling rate of 2048 the scale spaces are virtually identical.*

TABLE 1
*Comparison of scale spaces. For the "Blocks" data of [10] sampled in 64 points with a signal to noise ratio of 7, the eleven largest $\gamma$ values (see Lemma 4) for the deterministic signal (bottom) and the noisy signal (top) are compared. The last two values of the bottom row are equal to zero, since there are only nine ways to reconstruct the deterministic signal*

| 852 | 217 | 173 | 148 | 108 | 99.8 | 55.9 | 46.6 | 5.36 | 4.62 | 2.29 |
|-----|-----|-----|-----|-----|------|------|------|------|------|------|
| 885 | 249 | 159 | 142 | 100 | 99.1 | 80.2 | 41.3 | 38.9 | 0    | 0    |

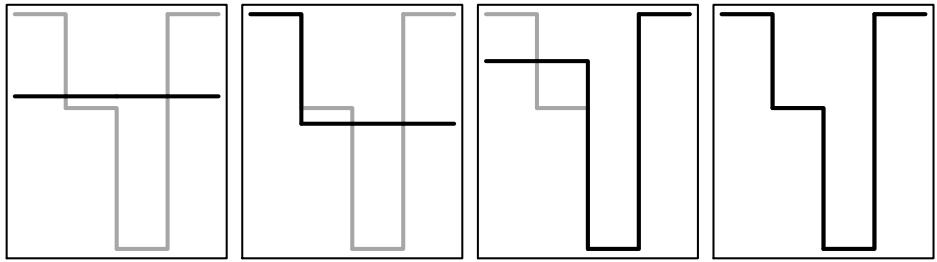

FIG 3. *Scale spaces of a sample function (grey line). The black lines show all reconstructions of the sample function for varying $\gamma$.*

**Discussion.** The scale space of a penalized estimator with hard thresholding type penalties generally does not have the same nice properties as its counterparts stemming from an $l_2$- or $l_1$-type penalty. In our case the function value at some point of the reconstruction does not change continuously or monotonically in the smoothing parameter. Moreover, the set of jumps of a best reconstruction with $k$ partitions is not necessarily contained in the set of jumps of a best reconstruction with $k'$ partitions for $k < k'$, see Figure 3. This leads to increased computational costs, as greedy algorithms in general do not yield an optimal solution. Indeed, one needs only $O(n \log n)$ steps to compute the estimate for a given $\gamma$ if the penalty is of $l_1$ type as in locally adaptive regression splines by [16], compared to $O(n^2)$ steps for the *Jplse*.

We mention, that penalizing the number of jumps corresponds to an $L^0$-penalty and is a limiting case of the [20] functional, when the dimension of the signal (image) is $d = 1$ [27], and results in "hard segmentation" of the data [24].

## 4. Proofs

**Some additional notation.** Throughout this section, we shorten $J(\hat{f}_n)$ to $J_n$. We set $S_n([0,1)) = \iota^n(\mathbb{R}^n)$, $\mathcal{B}_n = \sigma(S_n([0,1)))$. Observe that $\iota^n(\bar{f}_n)$ is just the conditional expectation $\mathbb{E}_{U_{0,1}}(f|\mathcal{B}_n)$, denoting the uniform distribution on $[0,1)$ by $U_{0,1}$. Similarly, for any finite $J \subset (0,1)$ define $\mathcal{B}_J = \sigma(\{[a,b) : a, b \in J \cup \{0,1\}\})$ and the partition $P_J = \{[a,b) : a, b \in J \cup \{0,1\}, (a,b) \cap J = \varnothing\}$. For our proofs it is convenient to formulate all minimization procedures on $L^2([0,1))$. Therefore we introduce the following functionals $\tilde{H}_\gamma^\infty$, $\tilde{H}_\gamma : L^2([0,1)) \times L^2([0,1)) \longmapsto \bar{\mathbb{R}}$, defined as

$$\tilde{H}_\gamma(g, f) = \begin{cases} \gamma \# J(g) + \|f - g\|^2 - \|f\|^2, & \text{if } g \in S_n([0,1)), \\ \infty, & \text{otherwise,} \end{cases}$$

$$\tilde{H}_\gamma^\infty(g, f) = H_\gamma^\infty(g, f) - \|f\|^2.$$

Clearly for each $f$, $\tilde{H}_\gamma^\infty$ has the same minimizers as $H_\gamma^\infty$, differing only by a constant. The following Lemma relates the minimizers of $\tilde{H}_\gamma$ and $H_\gamma$.

**Lemma 5.** *For all $f \in L^2([0,1))$ and $n \in \mathbb{N}$ we have $u \in \operatorname{argmin} H_\gamma(\cdot, \bar{f}_n)$ if and only if $\iota^n(u) \in \operatorname{argmin} \tilde{H}_\gamma(\cdot, f)$. Similarly, $u \in \operatorname{argmin} H_\gamma(\cdot, y)$ for $y \in \mathbb{R}^n$ if and only if $\iota^n(u) \in \operatorname{argmin} \tilde{H}_\gamma(\cdot, \iota^n(y))$.*

*Proof.* The second assertion follows from the fact that for $u, y \in \mathbb{R}^n$

$$\tilde{H}_\gamma(\iota^n(u), \iota^n(y)) = H_\gamma(u, y) - \|f\|^2.$$

Further, for $u \in \mathbb{R}^n$ we have $\langle \iota^n(\bar{f}_n) - f, \iota^n(\bar{f}_n) - \iota^n(u) \rangle = 0$ which gives

$$\begin{aligned}
\tilde{H}_\gamma(\iota^n(u), f) &= \gamma \# J(u) + \|f - g\|^2 - \|f\|^2 \\
&= \gamma \# J(u) + \|\iota^n(\bar{f}_n) - \iota^n(u)\|^2 + \|f - \iota^n(\bar{f}_n)\|^2 - \|f\|^2 \\
&= H_\gamma(u, \bar{f}_n) + \mathrm{const}_{f,n}
\end{aligned}$$

what completes the proof. □

The minimizers $g \in S([0,1))$ of $\tilde{H}_\gamma(\cdot, f)$ and $\tilde{H}_\gamma^\infty(\cdot, f)$ for $\gamma > 0$ are determined by their jump set $J(g)$ through the formula $g = \mathbb{E}_{U_{0,1}}(f|\mathcal{B}_{J(g)})$. In the sequel, we abbreviate

$$\mu_I(f) = \ell(I)^{-1} \int_I f(u)\, du$$

to denote the mean of $f$ on some interval $I$. In addition, we will use the abbreviation $f_J := \mathbb{E}_{U_{0,1}}(f|\mathcal{B}_J)$, such that for any partition $P_J$ of $[0,1)$

$$f_J = \sum_{I \in P_J} \mu_I(f) \mathbf{1}_I.$$

Further, we extend the noise in (1) to $L^2([0,1))$ by

$$\xi_n = \iota^n((\xi_{1,n}, \ldots, \xi_{n,n})).$$

### 4.1. Technical tools

We start by giving estimates on the behavior of $(\xi_n)_J = \sum_{I \in P_J} \mu_I(\xi_n) \mathbf{1}_I$.

**Lemma 6.** *Assume $(\xi_{i,n})_{n \in \mathbb{N}, 1 \leq i \leq n}$ satisfies (A). Then $\mathbb{P}$-almost surely for all intervals $I \subset [0,1)$ and all $n \in \mathbb{N}$*

$$\mu_I(\xi_n)^2 \leq \frac{\beta_n}{n\ell(I)}.$$

*Proof.* For intervals of the type $[(i-1)/n, j/n)$ with $i \leq j \in \mathbb{N}$ the claim is a direct consequence of (4). For general intervals, $[(i+p_1)/n, (j-p_2)/n)$ with $p_1, p_2 \in [0,1]$, we have to show that

$$(p_1 \cdot \xi_{i,n} + \xi_{i+1,n} + \cdots + \xi_{j-1,n} + p_2 \cdot \xi_{j,n})^2 - \beta_n(p_1 + p_2 + j - i - 1) \leq 0.$$

The left expression is convex over $[0,1]^2$ if it is considered as function in $(p_1, p_2)$. Hence it attains its maximum in an extreme point of $[0,1]^2$. □

**Lemma 7.** *There is a set of $\mathbb{P}$-probability one on which for all sequences $(J_n)_{n\in\mathbb{N}}$ of finite sets in $(0,1)$ the relation $\lim_{n\to\infty} \beta_n \# J_n/n = 0$ implies*

$$(\xi_n)_{J_n} \xrightarrow[n\to\infty]{L^2([0,1))} 0.$$

*Proof.* By Lemma 6 we find

$$(8) \qquad \|(\xi_n)_{J_n}\|^2 = \sum_{I\in P_{J_n}} \ell(I)\mu_I(\xi_n)^2 \leq \frac{\beta_n}{n}(\# J_n + 1),$$

This immediately gives the assertion. $\square$

Now we wish to show that the functionals epi-converge (see section 4.4). To this end we need two more results.

**Lemma 8.** *Let $(J_n)_{n\in\mathbb{N}}$ be a sequence of closed subsets in $(0,1)$ which satisfies the relation $\lim_{n\to\infty} \beta_n \# J_n/n = 0$. For $(g_n)_{n\in\mathbb{N}} \subset L^2([0,1))$ with $\|g_n - g\| \xrightarrow[n\to\infty]{} 0$, where $g_n$ is $\mathcal{B}_{J_n}$ measurable, we have almost surely*

$$\|f+\xi_n - g_n\|^2 - \|f+\xi_n\|^2 \xrightarrow[n\to\infty]{} \|f-g\|^2 - \|f\|^2.$$

*Proof.* First observe that

$$\|f+\xi_n - g_n\|^2 - \|f+\xi_n\|^2 = \|g_n\|^2 - 2\langle f, g_n\rangle - 2\langle \xi_n, g_n\rangle$$
$$= \|g_n\|^2 - 2\langle f, g_n\rangle - 2\langle (\xi_n)_{J_n}, g_n\rangle.$$

Since the sequence $(\|g_n\|)_{n\in\mathbb{N}}$ is bounded we can use Lemma 7 to deduce

$$\langle (\xi_n)_{J_n}, g_n\rangle \xrightarrow[n\to\infty]{\mathbb{P}\text{-a.s.}} 0.$$

This completes the proof. $\square$

Before stating the next result, we recall the definition of the Hausdorff metric $\rho_H$ on the space of closed subsets $CL(\Theta)$ of a compact metric space $(\Theta, \rho)$. For $\Theta' \subseteq \Theta \ni \vartheta$ we set

$$\operatorname{dist}(\vartheta, \Theta') = \inf\{\rho(\vartheta, \vartheta') : \vartheta' \in \Theta'\}.$$

Define

$$\rho_H(A,B) = \begin{cases} \max\{\sup_{x\in A}\operatorname{dist}(x,B), \sup_{y\in B}\operatorname{dist}(y,A)\}, & A, B \neq \varnothing, \\ 1, & A \neq B = \varnothing, \\ 0, & A = B = \varnothing, \end{cases}$$

With this metric, $CL(\Theta)$ is again compact for compact $\Theta$ [19, see].

**Lemma 9.** *The map*

$$L^2([0,1)) \ni g \mapsto \begin{cases} \# J(g), & g \in S([0,1)) \\ \infty, & g \notin S([0,1)) \end{cases} \in \mathbb{N} \cup \{0, \infty\}$$

*is lower semi-continuous, meaning the set $\{g \in S([0,1)) : \# J(g) \leq N\}$ is closed for all $N \in \mathbb{N} \cup \{0\}$.*

*Proof.* Suppose that $\|g_n - g\| \xrightarrow[n\to\infty]{} 0$ with $\#J(g_n) \leq N < \#J(g)$. Using compactness of the space of closed subsets $CL([0,1])$ and turning possibly to a subsequence, we could arrange that $J(g_n) \cup \{0,1\} \xrightarrow[n\to\infty]{} J \cup \{0,1\}$ for some closed $J \subset (0,1)$, where convergence is understood in Hausdorff metric $\rho_H$. Since the cardinality is lower semi-continuous with respect to the Hausdorff metric, $J$ must be finite. We conclude for $(s,t) \cap J = \varnothing$ and $\varepsilon > 0$ that $(s+\varepsilon, t-\varepsilon) \cap J(g_n) = \varnothing$ eventually, i.e. $g_n$ is constant on $(s+\varepsilon, t-\varepsilon)$. Next we observe that $g_n \mathbf{1}_{(s+\varepsilon,t-\varepsilon)}$ converges towards $g \mathbf{1}_{(s+\varepsilon,t-\varepsilon)}$ (in $L^2([0,1))$) what implies that $g$ is constant on $(s+\varepsilon, t-\varepsilon)$. Since $\varepsilon > 0$ was arbitrary, we derive that $g$ is constant on $(s,t)$. Consequently, $g$ is in $S([0,1))$ and $J(g) \subseteq J$. Using again lower semi-continuity of the cardinality in the space of compact subsets of $[0,1]$ shows that

$$\#J(g) > N \geq \limsup_n \#J(g_n) \geq \liminf_n \#J(g_n) \geq \#J \geq \#J(g).$$

This contradiction completes the proof. □

Now we can state the epi-convergence of $\tilde{H}_{\gamma_n}$ as function on $L^2([0,1))$.

**Lemma 10.** *For all sequences $(\gamma_n)_{n\in\mathbb{N}}$ satisfying (H) we have*

$$\tilde{H}_{\gamma_n}(\cdot, f + \xi_n) \xrightarrow[n\to\infty]{\text{epi}} \tilde{H}_\gamma^\infty(\cdot, f)$$

*almost surely. Here $\tilde{H}_{\gamma_n}, \tilde{H}_\gamma^\infty$ are considered as functionals on $L^2([0,1))$.*

*Proof.* We have to show that on a set with probability one we have

(i) If $g_n \xrightarrow[n\to\infty]{} g$ then $\liminf_{n\to\infty} \tilde{H}_{\gamma_n}(g_n, f + \xi_n) \geq \tilde{H}_\gamma^\infty(g, f)$.
(ii) For all $g \in L^2([0,1))$, there exists a sequence $(g_n)_{n\in\mathbb{N}} \subset L^2([0,1))$, $g_n \xrightarrow[n\to\infty]{} g$ with $\limsup_{n\to\infty} \tilde{H}_{\gamma_n}(g_n, f + \xi_n) \leq \tilde{H}_\gamma^\infty(g, f)$.

To this end, we fix the set where the assertions of Lemmas 7 and 8 hold simultaneously.

Ad 4.1: Without loss of generality, we may assume that $H_{\gamma_n}(g_n, f + \xi_n)$ converges in $\mathbb{R} \cup \infty$. If $g_n \notin S_n([0,1))$ for infinitely many $n$ or $\#J(g_n) > \tilde{H}_\gamma^\infty(g,f)/\gamma_n$ the relation 4.1 is trivially fulfilled. Otherwise, we obtain

$$\limsup_{n\to\infty} \frac{\beta_n}{n} \#J(g_n) \leq \limsup_{n\to\infty} \frac{\beta_n}{n\gamma_n} \tilde{H}_\gamma^\infty(g, f) = 0.$$

Hence we can apply Lemma 8. Together with Lemma 9 we obtain $\mathbb{P}$-a.s.

$$\liminf_{n\to\infty} \tilde{H}_{\gamma_n}(g_n, f + \xi_n)$$
$$\geq \liminf_{n\to\infty} \gamma_n J(g_n) + \liminf_{n\to\infty}(\|f + \xi_n - g_n\|^2 - \|f + \xi_n\|^2)$$
$$\geq \gamma J(g) + (\|f - g\|^2 - \|f\|^2) = \tilde{H}_\gamma^\infty(g, f).$$

Ad 4.1: If $g \notin S([0,1))$ and $\gamma > 0$ there is nothing to prove. If $\gamma = 0$ and still $g \notin S([0,1))$, choose $g_n$ as a best $L^2$-approximation of $g$ in $S_n([0,1))$ with at most $1/\sqrt{\gamma_n}$ jumps.

We claim that $\|g_n - g\| \to 0$ as $n \to \infty$. For that goal, let $\tilde{g}_{n,k}$ denote a best approximation of $g$ in $\{f \in S_n([0,1)) : \#J(f) \leq k\}$ and $\tilde{g}_k$ one in $\{f \in S([0,1)) : \#J(f) \leq k\}$.

Moreover, for every $n, k$ let $J_n^k \subset (0,1)$ be a perturbation of $J(\tilde{g}_k)$, with $nJ_n^k \in \mathbb{N}$, $\#J_n^k = \#J(\tilde{g}_k)$ and $\rho_H(J_n^k, J(\tilde{g}_k)) \leq 1/n$. Denote $g'_{n,k} = \tilde{g}_k \circ \lambda_{n,k}$ where $\lambda_{n,k} \in \Lambda_1$ fulfills $\lambda_{n,k}(J_n^k) = J(\tilde{g}_k)$. Since $(a,b) \mapsto \mathbf{1}_{[a,b)}$ is continuous in $L^2([0,1))$, we obtain readily $\|g'_{n,k} - \tilde{g}_k\| \to 0$. This implies for any $k \in \mathbb{N}$

$$\limsup_{n \to \infty} \|g_n - g\| \leq \limsup_{n \to \infty} \|\tilde{g}_{n,k} - g\| \leq \limsup_{n \to \infty} \|g'_{n,k} - g\| = \|\tilde{g}_k - g\|.$$

Since the right hand side can be made arbitrary small by choosing $k$, $g_n$ converges to $g$. Then Lemma 8 yields 4.1.

If $\gamma > 0$ and $g \in S([0,1))$, $g_n$ is chosen as a best approximation of $g$ in $S_n([0,1))$ with at most $\#J(g)$ jumps. Finally, in order to obtain 4.1, argue as before. $\square$

To deduce consistency with the help of epi-convergence, one needs to show that the minimizers are contained in a compact set. The following lemma will be applied to this end.

**Lemma 11.** *Assume $(\Theta, \rho)$ is a metric space. A subset $A \subset D([0,\infty), \Theta)$ is relatively compact if the following two conditions hold*

(B1) *For all $t \in \mathbb{R}_+$ there is a compact $K_t \subseteq \Theta$ such that*

$$g(t) \in K_t, \qquad (\text{for all } g \in A).$$

(B2) *For all $T > 0$ and all $\varepsilon > 0$ there exists a $\delta > 0$ such that for all $g \in A$ there is a step function $g_\varepsilon \in S([0,T), \Theta)$ such that*

$$\sup\{\rho(g(t), g_\varepsilon(t)) : t \in [0,T)\} < \varepsilon \qquad \text{and} \qquad \mathrm{mpl}(g_\varepsilon) \geq \delta,$$

*where mpl is the minimum distance between two jumps of $f \in S([0,T))$*

$$\mathrm{mpl}(f) := \min\{|s - t| : s \neq t \in J(f) \cup \{0, T\}\}.$$

*A subset $A \subset D([0,1), \Theta)$ is relative compact if the following two conditions hold*

(C1) *For all $t \in [0,1]$ there is a compact $K_t \subseteq \Theta$ such that*

$$g(t) \in K_t \qquad (\text{for all } g \in A).$$

(C2) *For all $\varepsilon > 0$ there exists a $\delta > 0$ such that for all $g \in A$ there is a step function $g_\varepsilon \in S([0,1), \Theta)$ such that*

$$\sup\{\rho(g(t), g_\varepsilon(t)) : t \in [0,1]\} < \varepsilon \qquad \text{and} \qquad \mathrm{mpl}(g_\varepsilon) \geq \delta.$$

*Proof.* We prove only the first assertion, as the proof of the second assertion can be carried out in the same manner.

According to [12], Theorem 6.3, it is enough to show that (B2) implies

$$\lim_{\delta \to 0} \sup_{g \in A} w_g(\delta, T) = 0$$

where

$$w_g(\delta, T) = \inf\left\{\max_{1 \leq i \leq v} \sup_{s,t \in [t_{i-1}, t_i)} \rho(g(s), g(t)) : \{t_1, \ldots, t_{v-1}\} \subset (0, T), \right.$$
$$\left. t_0 = 0, t_v = T, |t_i - t_j| > \delta\right\}.$$

So, fix $T > 0, \varepsilon > 0$ and choose $\delta$ from (B2). Then we set for $g \in A$ $\{t_0, \ldots, t_v\} = J(g_\varepsilon) \cup \{0, T\}$. Clearly, $\mathrm{mpl}(g_\varepsilon) > \delta$ implies $|t_i - t_j| > \delta$ for all $i \neq j$. For neighboring $t_{i-1}, t_i \in J(g_\varepsilon) \cup \{0, T\}$ and $s, t \in [t_{i-1}, t_i)$ we derive

$$\rho(g(s), g(t)) \leq \rho(g(s), g_\varepsilon(s)) + \rho(g_\varepsilon(s), g_\varepsilon(t)) + \rho(g_\varepsilon(t), g(t)) < \varepsilon + 0 + \varepsilon = 2\varepsilon.$$

This establishes the above condition and completes the proof. □

In the context of proving compactness we will also need the following result.

**Lemma 12.** *For any $f \in L^2([0,1))$ the set $\{f_J : J \subset (0,1), \#J < \infty\}$ is relatively compact in $L^2([0,1))$.*

*Proof.* The proof is done in several steps.
1. Since $(s,t) \mapsto \mathbf{1}_{[s,t)}$ is continuous,

$$\left\{ \sum_{i=1}^M \alpha_i \mathbf{1}_{I_i} : |\alpha_i| \leq z, I_i \subseteq [0,1) \text{ interval} \right\}$$

is the continuous image of a compact set and hence compact for all $M \in \mathbb{N}$ and $z > 0$.

2. If $f = \mathbf{1}_I$ for some interval $I$, we obtain for any $J \subset [0,1)$ that $f_J$ is a linear combination of at most three different indicator functions.

3. If $f = \sum_{i=1}^M \alpha_i \mathbf{1}_{I_i}$ is a step function and $J$ arbitrary then $f_J = \sum_{j=1}^{M'} \beta_j \mathbf{1}_{I'_j}$ holds by 2. for some $M' \leq 3M$. Using

$$\beta_j = \mu_{I'_j}(f) \leq \max_{i=1,\ldots,M} |\alpha_i|$$

as well as 1., we get that $\{f_J : J \subset [0,1)\}$ is relatively compact for step functions $f$.

4. Suppose $f \in L^2([0,1))$ is arbitrary and $\varepsilon > 0$. We want to show that we can cover $\{f_J : J \subset [0,1)\}$ by finitely many $\varepsilon$-balls. Fix a step function $g$ such that $\|f - g\| < \varepsilon/2$. By the Jensen Inequality for conditional expectations, we get $\|f_J - g_J\| < \varepsilon/2$ for all finite $J \subset [0,1)$. Further, by 3., there are finite sets $J_1, \ldots, J_p \subset [0,1)$ with $p < \infty$ such that $\min_{l=1,\ldots,p} \|g_J - g_{J_l}\| < \varepsilon/2$ for all finite $J \subset [0,1)$. This implies

$$\min_{l=1,\ldots,p} \|f_J - g_{J_i}\| \leq \min_{l=1,\ldots,p} \|g_J - g_{J_l}\| + \|f_J - g_J\| < \varepsilon$$

and the proof is complete. □

### 4.2. Behavior of the partial sum process

*Proof of Lemma 1.* The following Markov inequality is standard for triangular arrays fulfilling condition (A), [21], Section III, §4, and all numbers $\mu_i, i = 1, \ldots, n$:

$$\mathbb{P}(|\sum_{i=1}^n \mu_i \xi_{i,n}| \geq z) \leq 2 \exp\left(\frac{-z^2}{4\alpha n^\zeta \sum_i \mu_i^2}\right) \qquad \text{(for all } z \in \mathbb{R}\text{)}.$$

From this, we derive for $z^2 > 12\alpha$ that

$$\sum_{n \in \mathbb{N}} \sum_{1 \leq i \leq j \leq n} \mathbb{P}(|\xi_{i,n} + \cdots + \xi_{j,n}| \geq z\sqrt{j-i+1}\sqrt{n^\zeta \log n})$$

$$\leq 2 \sum_{n \in \mathbb{N}} n^2 e^{-\frac{z^2 \log n}{4\alpha}} = 2 \sum_n n^{-\frac{z^2 - 8\alpha}{4\alpha}} < \infty.$$

Hence, for $\epsilon > 0$ we have with probability one that

$$\max_{1 \le i \le j \le n} \frac{(\xi_{i,n} + \cdots + \xi_{j,n})^2}{(j - i + 1)} \ge (12 + \epsilon)\alpha n^\zeta \log n$$

only finitely often. □

For the proof of Lemma 2, we need an auxiliary lemma. Denote by

$$\mathcal{D}_n = \big\{(i,j) : 1 \le i \le j \le n \quad \text{such that} \quad i = k2^l, j = (k+1)2^l \\ \text{for some} \quad l, k \in \{0, 1, 2, \ldots\} \big\}$$

the set of all pairs $(i, j)$ which are endpoints of dyadic intervals contained in $\{1, \ldots, n\}$.

**Lemma 13.** *Assume $x \in \mathbb{R}^n$ such that*

(9) $$\max_{(i,j) \in \mathcal{D}_n} \frac{|x_i + \cdots + x_j|}{\sqrt{j - i + 1}} \le c$$

*for some $c > 0$. Then*

$$\max_{1 \le i \le j \le n} \frac{|x_i + \cdots + x_j|}{\sqrt{j - i + 1}} \le (2 + \sqrt{2})c.$$

*Proof.* Without loss of generality we may assume that $n = 2^m$ for some $m \in \mathbb{N}$ (and add some zeros otherwise). First, we prove by induction on $m$ that (9) implies

(10) $$\max_{1 \le j \le n} \frac{|x_1 + \cdots + x_j|}{\sqrt{j}} \le (1 + \sqrt{2})c.$$

For $m = 0$ there is nothing to prove. Now assume that the statement is true for $m$. Let $2^m < j \le 2^{m+1}$. Note that

$$\frac{|x_1 + \cdots + x_j|}{\sqrt{j}} \le \frac{\sqrt{2^m}}{\sqrt{j}} \frac{|x_1 + \cdots + x_{2^m}|}{\sqrt{2^m}} + \frac{\sqrt{j - 2^m}}{\sqrt{j}} \frac{|x_{2^m+1} + \cdots + x_j|}{\sqrt{j - 2^m}}.$$

Apply the induction hypothesis to the second summand to obtain

$$\frac{|x_1 + \cdots + x_j|}{\sqrt{j}} \le \frac{\left(\sqrt{2^m} + (1 + \sqrt{2})\sqrt{j - 2^m}\right)}{\sqrt{j}} c.$$

For $2^m + 1 \le j \le 2^{m+1}$ the expression on the right hand side is maximal for $j = 2^{m+1}$ with maximum $(1 + \sqrt{2})c$. Hence the statement holds also for $m + 1$ and we have shown that (9) implies (10).

The claim is again proven by induction on $m$. For $m = 0$ there is nothing to prove. Assume that the statement is true for $m$. If $i \le j \le 2^m$ or $2^m < i \le j \le 2^{m+1}$ the statement follows by application of the induction hypotheses to $(x_1, \ldots, x_{2^m})$ and $(x_{2^m+1}, \ldots, x_{2^{m+1}})$, respectively. Now suppose $i < 2^m < j$. Then

$$\frac{|x_i + \cdots + x_j|}{\sqrt{j - i + 1}} \le \frac{\sqrt{2^m - i + 1}}{\sqrt{j - i + 1}} \frac{|x_i + \cdots + x_{2^m}|}{\sqrt{2^m - i + 1}} \\ + \frac{\sqrt{j - 2^m}}{\sqrt{j - i + 1}} \frac{|x_{2^m+1} + \cdots + x_j|}{\sqrt{j - 2^m}}$$

Application of (10) to $x' = (x_{2^m}, x_{2^m-1}, \ldots, x_1)$ and $\tilde{x} = (x_{2^m+1}, \ldots, x_{2^{m+1}})$ then gives

$$\frac{|x_i + \cdots + x_j|}{\sqrt{j-i+1}} \leq \frac{\sqrt{2^m - i + 1} + \sqrt{j - 2^m}}{\sqrt{j-i+1}} \ (1 + \sqrt{2})c \leq \sqrt{2}(1 + \sqrt{2})c. \quad \square$$

*Proof of Lemma 2.* [8] show that for $m \geq 1$ and some constant $C_m$ depending on $m$ only

$$\mathbb{E}\Big(\frac{|\xi_{i,n} + \cdots + \xi_{j,n}|^{2m}}{(j-i+1)^m}\Big) \leq C_m \frac{\mathbb{E}|\xi_{i,n}|^{2m} + \cdots + \mathbb{E}|\xi_{j,n}|^{2m}}{j-i+1}.$$

The Markov inequality then yields for any $z > 0$ and all $1 \leq i \leq j \leq n$

$$\mathbb{P}\Big(\frac{|\xi_{i,n} + \cdots + \xi_{j,n}|}{\sqrt{j-i+1}} \geq z\Big) \leq \frac{C_m \sup_{i,n} \mathbb{E}|\xi_{i,n}|^{2m}}{z^{2m}}.$$

Since there are at most $2n$ dyadic intervals contained in $\{1, \ldots, n\}$, we obtain by Lemma 13 for any $C > 0$ that

$$\sum_{n \in \mathbb{N}} \sum_{1 \leq i \leq j \leq n} \mathbb{P}\Big(\frac{|\xi_{i,n} + \cdots + \xi_{j,n}|}{\sqrt{j-i+1}} \geq C(n \log n)^{1/m}\Big)$$

$$\leq \sum_{n \in \mathbb{N}} \sum_{(i,j) \in \mathcal{D}_n} \mathbb{P}\Big(\frac{|\xi_{i,n} + \cdots + \xi_{j,n}|}{\sqrt{j-i+1}} \geq (2 + \sqrt{2})C(n \log n)^{1/m}\Big)$$

$$\leq \frac{C_m \sup_{i,n} \mathbb{E}|\xi_{i,n}|^{2m}}{(2 + \sqrt{2})^{2m} C^{2m}} \sum_{n \in \mathbb{N}} \frac{2n}{n^2 \log^2 n} < \infty.$$

The claim follows by application of the Borel–Cantelli lemma. $\quad \square$

### 4.3. Consistency of the estimator

The proofs in this section use the concept of epi-convergence. It is introduced in Appendix.

*Proof of Lemma 3.* For $\gamma = 0$ there is nothing to prove. Assume $\gamma > 0$ and $g \in S([0,1))$ with $\#J(g) > \|f\|^2/\gamma$. This yields

$$H_\gamma^\infty(0, f) = \|f\|^2 < H_\gamma^\infty(g, f).$$

Moreover, observe that for $g \in S([0,1))$ we have $H_\gamma^\infty(g, f) \geq H_\gamma^\infty(f_{J(g)}, f)$. Thus, it is enough to regard the set $\{f_J : \#J \leq \|f\|^2/\gamma\}$, which is relatively compact in $L^2([0,1))$ by Lemma 12. This proves the existence of a minimizer. $\quad \square$

*Proof of Theorem 1 and Theorem 2.* By the reformulation of the minimizers in Lemma 5, Lemma 10 and Theorem 5 (see Appendix) it is enough to prove that almost surely there is a compact set containing

$$\bigcup_{n \in \mathbb{N}} \operatorname{argmin} \tilde{H}_{\gamma_n}(\cdot, f + \xi_n).$$

First note that all $f_n \in \operatorname{argmin} \tilde{H}_{\gamma_n}(\cdot, f + \xi_n)$ have the form $(f + \xi_n)_{J_n}$ for some (random) sets $J_n$. Comparing $\tilde{H}_{\gamma_n}(f_n, f + \xi_n)$ with $\tilde{H}_{\gamma_n}(0, f + \xi_n) = 0$, we obtain the a priori estimate

$$\gamma_n \# J_n \leq \|(f + \xi_n)_{J_n}\|^2 \leq 2\|f\|^2 + 2\,\|(\xi_n)_{J_n}\|^2 \leq 2\|f\|^2 + \frac{2\beta_n}{n}(\#J_n + 1)$$

for all $n \in \mathbb{N}$. Since $\gamma_n > \frac{4\beta_n}{n}$ eventually, we find $\mathbb{P}$-a.s.

$$\#J_n \leq \frac{2\|f\|^2 + \frac{2\beta_n}{n}}{\gamma_n - \frac{2\beta_n}{n}} = \mathrm{O}(\gamma_n^{-1}).$$

Application of Lemma 7 gives $\lim_{n\to\infty} (\xi_n)_{J_n} = 0$ almost surely. Since by Lemma 12, $\{f_{J_n} : n \in \mathbb{N}\}$ is relatively compact in $L^2([0,1))$, relative compactness of the set $\bigcup_{n\in\mathbb{N}} \operatorname{argmin} \tilde{H}_{\gamma_n}(\cdot, f + \xi_n)$ follows immediately. This completes the proofs. □

*Proof of Theorem 3, part* (i) . Theorem 1 and Lemma 9 imply

$$\liminf_{n\to\infty} \#J_n \geq \#J(f_\gamma).$$

Suppose $\limsup_{n\to\infty} \#J_n \geq \#J(f_\gamma) + 1$. Let $f_{\gamma,n}$ be an approximation of $f_\gamma$ from $S_n([0,1))$ with the same number of jumps as $f_\gamma$. Then we could arrange $f_{\gamma,n} \xrightarrow[n\to\infty]{} f_\gamma$ such that $\lim_{n\to\infty} \tilde{H}_\gamma(f_{\gamma,n}, f + \xi_n) = \tilde{H}_\gamma^\infty(f_\gamma, f)$. Moreover, we know

$$\limsup_{n\to\infty} \tilde{H}_\gamma(\hat{f}_n, f + \xi_n) \geq \gamma + \tilde{H}_\gamma^\infty(f_\gamma, f) = \gamma + \lim_{n\to\infty} \tilde{H}_\gamma(f_{\gamma,n}, f + \xi_n)$$

which contradicts that $\hat{f}_n$ is a minimizer of $\tilde{H}_\gamma(\cdot, f + \xi_n)$ for all $n$. Therefore, $\#J_n = \#J(f_\gamma)$ eventually.

Next, chose by compactness a subsequence such that $J_n \cup \{0,1\}$ converges in $\rho_H$. Then, by Lemma 9, the limit must be $J(f_\gamma) \cup \{0,1\}$. Consequently, the whole sequence $(J_n)_{n\in\mathbb{N}}$ converges to $J(f_\gamma)$ in the Hausdorff metric.

Thus eventually, there is a 1-1 correspondence between $P_{J_n}$ and $P_{J(f_\gamma)}$ such that for each $[s,t] \in P_{J(f_\gamma)}$ there are $[s_n, t_n] \in P_{J_n}$ with

$$s_n \xrightarrow[n\to\infty]{} s \qquad \text{and} \qquad t_n \xrightarrow[n\to\infty]{} t.$$

By Lemma 6 and continuity of $(s,t) \mapsto \mathbf{1}_{[s,t)}$, we find

$$\mu_{[s_n, t_n)}(f + \xi_n) \xrightarrow[n\to\infty]{} \mu_{[s,t)}(f).$$

Construct $\lambda_n \in \Lambda_1$ linearly interpolating $\lambda_n(s_n) = s$. Then

$$\mathrm{L}(\lambda_n) \xrightarrow[n\to\infty]{} 1$$

as well as

$$\|\hat{f}_n - f \circ \lambda_n\|_\infty = \max_{I \in P_{J(f_\gamma)}} |\mu_{\lambda^{-1}(I)}(f + \xi_n) - \mu_I(f)| \xrightarrow[n\to\infty]{} 0$$

which completes the proof. □

*Proof of Theorem 3, part* (ii). The proof can be carried out in the same manner as the proof of Theorem 4, part (ii) in [2]. The only difference is, that it is necessary to attend the slightly different rates of the partial sum process (4). □

## *4.4. Convergence of scale spaces*

*Proof of Lemma 4.* It is clear, that each $g \in \operatorname{argmin} H_\gamma^\infty(\cdot, f)$ is determined by its jump set. Further, if $g_1, g_2 \in S([0,1))$ with $\#J(g_1) = \#J(g_2)$ and $\|f-g_1\| = \|f-g_2\|$ then $g_1$ is a minimizer of $H_\gamma^\infty(\cdot, f)$ if and only if $g_2$ is.

Since $H_\gamma^\infty(0, f) = \|f\|^2$ we have that $\gamma \in [\nu, \infty)$ implies $J(g) \leq \|f\|^2/\nu$, for a minimizer $g$ of $H_\gamma^\infty(\cdot, f)$. Hence on $[\nu, \infty)$ we have that

$$\min H_\gamma^\infty(\cdot, f) = \min\{k\gamma + \Delta_k(f) : k \leq \|f\|^2/\nu\}$$

with $\Delta_k(f)$ defined by

$$\Delta_k(f) := \inf\{\|g - f\| : g \in S([0,1)), \#J(g) \leq k\}.$$

For each $\nu$ the map $\gamma \mapsto \min H_\gamma^\infty(\cdot, f)$ is thus a minimum of a finite collection of linear functions with pairwise different slopes on $[\nu, \infty)$. If there are different $k, k'$ and $\gamma$ with $k\gamma + h_k = k'\gamma + h_{k'}$ it follows $\gamma = (h_{k'} - h_k)/(k - k')$. From this it follows that there are only finitely many $\gamma$ where $\#\{k : k\lambda + \Delta_k(f) = \min H_\gamma^\infty(\cdot, f)\} > 1$. Further, $\operatorname{argmin} H_\gamma^\infty(\cdot, y)$ is completely determined by the $k$ which realize this minimum. Call those $\gamma$, for which different $k$ realize the minimum, changepoints of $\gamma \mapsto \min H_\gamma^\infty(\cdot, f)$. Since the above holds true for each $\nu > 0$, there are only countably many changepoints in $[0, \infty)$. This completes the proof. □

*Proof of Theorem 4.* It is easy to see that the assumptions imply $J(\tau) = \{\gamma_m : m = 1, \ldots, m(x)\}$ for the sequence $(\gamma_m)_{m=0}^{m(x)} \subset \mathbb{R} \cup \infty$ of Lemma 4. Since the scale space $\tau$ is uniquely determined by its jump points, this proves the uniqueness claim.

For the proof of the almost sure convergence, note that Theorem 1 and Theorem 3, part (i) show that $\hat\tau_n(\zeta) \to_{n \to \infty} \tau(\zeta)$ if $\zeta$ is a point of continuity of $\tau$, i.e. $\#\operatorname{argmin} H_{1/\zeta}^\infty(\cdot, f) = 1$. Convergence in all continuity points together with relative compactness of the sequence implies convergence in the Skorokhod topology. Hence, it is enough to show that $\{\hat\tau_n : n \in \mathbb{N}\}$ is relatively compact.

To this end, we will use Lemma 11. In the proof of Theorem 1 it was shown, that the sequence $(T_\gamma(Y_n))_{n \in \mathbb{N}}$ is relatively compact in $L^2(0,1)$. To prove relative compactness in $D([0,1))$ we follow the lines of the proof of Theorem 3, part (i). Similarly we find that

$$\limsup_{n \to \infty} \#J_n \leq \max_{g \in \operatorname{argmin} H_{1/\zeta}^\infty(\cdot, f)} \#J(g).$$

For each subsequence of $(T_\gamma(Y_n))_{n \in \mathbb{N}}$, consider the subsequence of corresponding jump sets. By compactness of $CL([0,1])$ we choose a converging sub-subsequence and argue as in the proof mentioned above that the corresponding minimizers converge to a limit in $\operatorname{argmin} H_{1/\zeta}^\infty(\cdot, f)$. Thus we have verified condition (B1).

For the proof of (B2), we will show by contradiction that for all $T > 0$ we have

$$\inf\{\operatorname{mpl}(\hat\tau_n|_{[0,T]}) : n \in \mathbb{N}\} > 0.$$

This, obviously, would imply (B2). Observe that $\hat\tau_n$ jumps in $\zeta$ only if there are two jump sets $J \neq J'$ such that $H_{1/\zeta}((Y_n)_J, Y_n) = H_{1/\zeta}((Y_n)_{J'}, Y_n)$ and $H_{1/\zeta}((Y_n)_J, Y_n) \leq H_{1/\zeta}((Y_n)_{J''}, Y_n)$ for all $J''$.

If (B2) is not fulfilled for $(\tau_n)_{n \in \mathbb{N}}$, we can switch by compactness to a subsequence and find sequences $(\zeta^1_n)_{n \in \mathbb{N}}, (\zeta^2_n)_{n \in \mathbb{N}}$ with $\zeta^1_n, \zeta^2_n \in J(\hat\tau_n)$, $\zeta^1_n < \zeta^2_n$ and $\zeta^1_n \xrightarrow[n \to \infty]{} \zeta$, $\zeta^2_n \xrightarrow[n \to \infty]{} \zeta$ for some $\zeta \geq 0$. Choosing again a subsequence, we could assume that the jump sets $J_n^1, J_n^2, J_n^3$ of minimizers $\hat{f}_n^k \in \iota^n(\operatorname{argmin} H_{\gamma_n}(\cdot, Y_n))$ for some sequences $\gamma_n^1 - 1/\zeta_n^1 \downarrow 0$, $\gamma_n^2 \in (1/\zeta_n^2, 1/\zeta_n^1)$ and $\gamma_n^3 - 1/\zeta_n^2 \uparrow 0$ are constant and $(\hat{f}_n^k)_{n \in \mathbb{N}}$, $k = 1, 2, 3$, converge. Further, we know from this choice of $\gamma_n^k$ and Lemma 4 that $\#J_n^1 > \#J_n^2 > \#J_n^3$. This implies

$$(11) \qquad \gamma_n^1 + \gamma_n^2 + \|\iota^n(Y_n) - \hat{f}_n^1\|^2 < \gamma_n^2 + \|\iota^n(Y_n) - \hat{f}_n^2\|^2 < \|\iota^n(Y_n) - \hat{f}_n^3\|^2.$$

The same arguments as in Theorem 1 and Theorem 3, part (i) respectively, yield $\{\lim_{n \to \infty} \hat{f}_n^k : k = 1, 2, 3\} \subseteq \operatorname{argmin} H_{1/\zeta}^\infty(\cdot, f)$. Since (11) holds for all $n$, the limits are pairwise different. This contradicts $\# \operatorname{argmin} H_{1/\zeta}^\infty(\cdot, x) \leq 2$ and proves (B2).

Thus $\{\hat{\tau}_n : n \in \mathbb{N}\}$ is relatively compact in $D([0, \infty), D([0, 1)))$ as well as in $D([0, \infty), L^2[0, 1])$ and the proof is complete. □

**Appendix: Epi-Convergence**

Instead of standard techniques from penalized maximum likelihood regression, we use the concept of epi-convergence (see for example [5, 13]). This allows for simple formulation and more structured proofs. The main arguments to derive consistency of estimates which are (approximate) minimizers for a sequence of functionals can briefly be summarized by epi-convergence + compactness + uniqueness a.s. $\Rightarrow$ strong consistency.

We give here the definition of epi-(or $\Gamma$-)convergence together with the results from variational analysis which are relevant for the subsequent proofs.

**Definition 1.** Let $F_n : \Theta \longmapsto \mathbb{R} \cup \infty$, $n = 1, \ldots, \infty$ be numerical functions on a metric space $(\Theta, \rho)$. $(F_n)_{n \in \mathbb{N}}$ *epi-converges to* $F_\infty$ *(Symbol $F_n \xrightarrow[n \to \infty]{\text{epi}} F_\infty$)* if

(i) for all $\vartheta \in \Theta$, and sequences $(\vartheta_n)_{n \in \mathbb{N}}$ with $\vartheta_n \xrightarrow[n \to \infty]{} \vartheta$

$$F_\infty(\vartheta) \leq \liminf_{n \to \infty} F_n(\vartheta_n)$$

(ii) for all $\vartheta \in \Theta$ there exists a sequence $(\vartheta_n)_{n \in \mathbb{N}}$ with $\vartheta_n \xrightarrow[n \to \infty]{} \vartheta$ such that

$$(12) \qquad F_\infty(\vartheta) \geq \limsup_{n \to \infty} F_n(\vartheta_n)$$

The main, useful conclusions from epi-convergence are given by the following theorem.

**Theorem 5** ([1], Theorem 5.3.6). *Suppose $F_n \xrightarrow[n \to \infty]{\text{epi}} F_\infty$.*

(i) *For any converging sequence $(\vartheta_n)_{n \in \mathbb{N}}$, $\vartheta_n \in \operatorname{argmin} F_n$, it holds necessarily $\lim_{n \to \infty} \vartheta_n \in \operatorname{argmin} F_\infty$.*
(ii) *If there is a compact set $K \subset \Theta$ such that $\varnothing \neq \operatorname{argmin} F_n \subset K$ for large enough $n$ then $\operatorname{argmin} F_\infty \neq \varnothing$ and*

$$\operatorname{dist}(\vartheta_n, \operatorname{argmin} F_\infty) \xrightarrow[n \to \infty]{} 0$$

*for any sequence $(\vartheta_n)_{n \in \mathbb{N}}$, $\vartheta_n \in \operatorname{argmin} F_n$.*
(iii) *If, additionally, $\operatorname{argmin} F_\infty$ is a singleton $\{\vartheta\}$ then*

$$\vartheta_n \xrightarrow[n \to \infty]{} \vartheta$$

*for any sequence $(\vartheta_n)_{n \in \mathbb{N}}$, $\vartheta_n \in \operatorname{argmin} F_n$.*